\documentclass[a4paper,11pt]{article}
\usepackage{latexsym}
\usepackage{amsfonts,amssymb}
\usepackage{amscd}
\usepackage[dvips]{graphicx}

\font\erm=cmr9

\author{M.~Dziemia{\'n}czuk}
\title{On Cobweb admissible sequences}

\newtheorem{defn}{Definition}
\newtheorem{problem}{Problem}
\newtheorem{theoremn}{Theorem}
\newtheorem{observen}{Observation}
\newtheorem{algorithm}{Algorithm}

\newtheorem{lemma}{Lemma}

\begin{document}
\begin{center}
\noindent {\Large \textsc{On Cobweb Admissible Sequences}}  \\ 
			\textit{ The Production Theorem}

\vspace{0.5cm}

\noindent {Maciej Dziemia\'nczuk}

\vspace{0.5cm}

\noindent {\erm Student in the Institute of Computer Science, Bia\l ystok University (*)}

\noindent {\erm PL-15-887 Bia\l ystok, st. Sosnowa 64, Poland}

\noindent {\erm e-mail: Maciek.Ciupa@gmail.com}

\noindent {\erm (*) former Warsaw University Division}
\end{center}

\vspace{1cm}

\noindent \textbf{Summary} 

\noindent In this note further  clue decisive  observations on cobweb admissible sequences are shared with the audience. In particular  an announced  proof of the Theorem 1 (by Dziemia\'nczuk)  from \cite{1}  announced in India -Kolkata- December 2007 is delivered here.  Namely here and there we claim that any cobweb admissible sequence F is at the point product of primary cobweb admissible sequences taking values one and/or certain power of an appropriate primary number $p$.

\noindent Here also an algorithm  to produce the family of all cobweb-admissible sequences i.e. the Problem 1 from  \cite{1} i.e. one of several problems posed in source papers \cite{2,3} is solved  using the idea and methods implicitly present  already in \cite{4}.

\vspace{0.4cm}
\noindent  Presented at Gian-Carlo Polish Seminar:

\noindent \emph{http://ii.uwb.edu.pl/akk/sem/sem\_rota.htm}

\section{Preliminaries}

\noindent The notation from \cite{2,3,1} is being here taken for granted.

\begin{defn}[\cite{2,3,1}] \label{def:admissible}
A sequence $F$ is called cobweb-admissible iff for any $n,k\in \mathbb{N} \cup\{0\}$

\begin{equation}
	{n \choose k}_F = \frac{n_F\cdot(n-1)_F\cdot...\cdot(n-k+1)_F}{1_F\cdot 2_F\cdot...\cdot k_F} \in \mathbb{N} 
\end{equation}
\end{defn}

\vspace{0.2cm}

\begin{problem}[\cite{2,3,1}] 
	Find effective characterizations and/or an algorithm to produce the cobweb admissible sequences i.e. find all examples.
\end{problem}

\section{Primary cobweb admissible sequence}

\noindent Throughout this  paper we shall consequently use $p$ letter only for primary numbers.

\vspace{0.4cm}
\begin{defn} \label{def:primary}
	A cobweb admissible sequence $P(p)\equiv\{n_P\}_{n\geq 0}$ valued one and/or powers of one certain primary number $p$ i.e. $n_P\in\{1,p,p^2,p^3,...\}$ is called primary cobweb admissible sequence.
\end{defn}

\begin{theoremn}\label{th:1}
	Any cobweb-admissible sequence $F$ is at the point product of primary cobweb-admissible sequences $P(p)$.
\end{theoremn}

\noindent \textbf{P}\textbf{\footnotesize{ROOF}}

\noindent Given any cobweb admissible sequence $F=\{n_F\}_{n\geq 0}$,  each of its elements can be represented as a product of primary numbers' powers i.e. $\break n_F = \prod_{s\geq 1}{p_s^{\alpha(n,s)}}$. 
Therefore the sequence $F$ is at the point product of sequences $P(p_1), P(p_2), ...$ such that  $P(p_s)\equiv\{n_{P_s}\}_{n\geq 0}$ and $n_{P_s} = p_s^{\alpha(n,s)}$. Each of primary sequences $P(p_s)$ where $s=1,2,3,...$ is cobweb admissible as following holds for any $n,k\in \mathbb{N} \cup\{0\}$

\begin{equation}
	\begin{array}{l}
	{n \choose k}_{P(p_1)} \cdot {n \choose k}_{P(p_2)} \cdot ... = 
	{n \choose k}_{F} \in \mathbb{N}  \Rightarrow
	\\ 
	\Rightarrow 
	{n \choose k}_{P(p_s)} = \frac{n_{P(p_s)}^{\underline{k}}}{k_{P(p_s)}!} = 
	\frac{p_s^N}{p_s^K} \in \mathbb{N} 
	\end{array}
\end{equation}

\noindent where $N$ stands for the sum of index powers' of primary numbers $p_s$  in  $n_F,\nobreak{(n-1)_F},...,(n-k+1)_F$ product expansion via primary numbers and correspondingly $K$ is  the  index powers' sum for $k$ first elements of the sequence $F$ $\blacksquare$

\section{Primary cobweb admissible sequences family}

In this section we define a family $\mathcal{A}(p)$ of all primary cobweb admissible sequences taking values one and/or certain power of an appropriate primary number $p$ . In the next part of this section we present the family in the graph structure of a tree defined in algorithmic way in what follows.
\vspace{0.4cm}

\noindent For this aim let consider a primary cobweb admissible sequence $F \equiv P(p)$ and its corresponding family of sequences  $B(F)\equiv\{n_{B(F)}\}_{n\geq 0}$ such that 

$$
	n_{B(F)} = m \leftrightarrow n_F = p^m
$$

\vspace{0.2cm}
\noindent In the sequel we shall consider sequences for arbitrary but fixed  one primary number $p$, therefore we use abbreviatio $P(p)\equiv P$.

\begin{lemma}\label{lem:1}
Natural number valued sequence $F\equiv\{n_F\}_{n\geq 0}$ is primary cobweb admissible $P(p)$ iff for any natural number $n$, $n_F\in\{1,p,p^2,p^3,...\}$ and
$$
	\forall_{1\leq k\leq \lfloor n/2 \rfloor} \sum_{s=n-k+1}^{n}{s_{B(F)}} \geq \sum_{s=1}^{k}{s_{B(F)}}.
$$
\end{lemma}

\noindent \textbf{P}\textbf{\footnotesize{ROOF}}

\noindent \textbf{The first steep}. Given any primary cobweb admissible sequence $\break F\equiv\{n_F\}_{n\geq 0}$. From the Definition \ref{def:primary} we know that $n_F \in \{1,p,p^2,...\}$ for certain primary number $p$. From the Definition \ref{def:admissible} we readily infer that for any $n,k\in \mathbb{N} \cup\{0\}$
$$
	{n \choose k}_F = 
	\frac{p^{n_{B(F)}}\cdot...\cdot p^{(n-k+1)_{B(F)}}}{p^{1_{B(F)}} \cdot p^{2_{B(F)}}\cdot...\cdot p^{k_{B(F)}}} = 
	\frac{p^N}{p^K} \in \mathbb{N} \Rightarrow N \geq K
$$

\noindent where $N =  \sum_{s=n-k+1}^{n}{s_{B(F)}}$ and $K = \sum_{s=1}^{k}{s_{B(F)}}$.

\vspace{0.4cm}
\noindent \textbf{The second steep}. Given any sequence $F\equiv\{n_F\}_{n\geq 0}$ where $\break{n_F\in\{1,p,p^2,p^3,...\}}$ and $\forall_{1\leq k\leq \lfloor n/2 \rfloor} \sum_{s=n-k+1}^{n}{s_{B(F)}} \geq \sum_{s=1}^{k}{s_{B(F)}}$ (*). Then for any natural $n,k$ below takes place
$$
	{n \choose k}_F = \frac{p^{n_{B(F)}}\cdot...\cdot p^{(n-k+1)_{B(F)}}}{p^{1_{B(F)}} \cdot p^{2_{B(F)}}\cdot...\cdot p^{k_{B(F)}}} = C \wedge (*) \Rightarrow C \in \mathbb{N} 
$$
\noindent $\blacksquare$

\begin{defn}[Primary cobweb admissible tree]\label{def:tree}
Let $G(p)$ to be a \break weighted tree $G(p) = \langle V,E,\delta \rangle$ where $V$ stays for set of vertices, $E$ denotes a set of nodes and function $\delta$ which assigns weight for any vertex $v\in V$ such that $\delta(v)\in\{0,1,2,...\}$ and $p$ - primary number. We shall define the corresponding graph $G(p)$ via the following recurrence:
\end{defn}

\begin{enumerate}
	\item {$v_0\in V$ is called the root with weight $\delta(v_0)=0$}
	\item {If $(v_0,v_1,...,v_{n-1})$ is a path of graph $G(p)$ then $(v_0,v_1,...,v_{n-1},v_{n})$ is too if, and only if $\forall_{1\leq k \leq \lfloor n/ 2\rfloor}{N_{n,k} \geq K_k} )$}
\end{enumerate}

\noindent where $N_{n,k} = \sum_{i=n-k+1}^{n}{\delta(v_{i})}$ and $K_k = \sum_{i=1}^{k}{\delta(v_{i})}$ .

\vspace{0.4cm}

\noindent \textbf{Conclusion 1}\\
Any path $(v_0,v_1,...,v_n)$ from the root $v_0$ to vertex $v_n$ encodes the first $n$ terms with $0_{F}$ of primary cobweb admissible sequence $F \equiv \{n_F\}_{\geq 0}$, $\nobreak{n_F \in \{1,p,p^2,...\}}$ with help of elements' exponent powers' sequence $B(P)$ such that $k_{B(F)} = \delta(v_k)$ i.e. the $n+1$-tuple 
$$
(v_0,v_1,...,v_n) \in V^{n+1} \leftrightarrow (\delta(v_0),\delta(v_1),...,\delta(v_n))
$$
exactly encodes finite primary cobweb admissible sequence $F$ valued by one and/or powers of primary number $p$.

\begin{observen}\label{obs:nonempty}
	If any path $(v_0,v_1,...,v_n)$ encodes $n$ the first terms with $0_F$ of primary cobweb admissible sequence $F$ then there exists infinite number of successors vertices $v_{n+1}$ which encode primary cobweb admissible sequence $F'$ specified by these $n$ first terms with $0_F$ and the one additional $(n+1)_{F'} = \delta(v_{n+1})$ term.
\end{observen}

\noindent \textbf{P}\textbf{\footnotesize{ROOF}}

\noindent If any path $(v_0,v_1,...,v_n)$ encodes $n$ the first terms with $0_F$ of primary cobweb admissible sequence $F$ then there exists infinite number of natural numbers $M$ such that $\delta(v_{n+1}) = M$ and $N_{n,k-1}+M \geq K_k$.

\vspace{0.4cm}
\noindent Consequently, now we present an algorithm to generate primary cobweb admissible tree.

\begin{algorithm}[primary cobweb-admissible tree] \label{alg:tree}
We shall begin with the root $v_0$ of graph $G(p)$ from Definition \ref{def:tree} and in the next steeps, from any path $(v_0,v_1,...,v_n)$ we obtain the very next one  $(v_0,v_1,...,v_n,v_{n+1})$.
\end{algorithm}

\noindent \textbf{Input:} Any path $(v_0,v_1,...,v_{n})$ of $G(p)$ which encodes $n$ the first terms with $0_F$ of primary cobweb admissible sequence $F$.

\vspace{0.4cm}

\noindent \textbf{Output:} Non-empty set $\emptyset \neq \Delta_{n} \subseteq \{v_{n+1}:\delta(v_{n+1}) \in \{0,1,2,...\}\}$ with vertices' successors for $v_n$ vertex such that the paths $(v_0,v_1,...,v_n,v_{n+1})$ where $\nobreak{v_{n+1}\in\Delta_{n+1}}$  encodes primary cobweb-admissible sequence, too.

\vspace{0.4cm}
\noindent Under the convenient notation for vertices $v(s)\equiv v_{n+1} \wedge \delta(v_{n+1}) = s$ note now the following. 

\vspace{0.4cm}
\noindent \textbf{Steps:}

\begin{enumerate}
	\item {If $n = 1 \\
		\Delta_1 = \{v(0), v(1), v(2), ...\}$}
	\item {If $n = 2 \\
		\Delta_2 = \{ v(m), v(m+1), v(m+2), ... \}$, where  $m=\delta(v_1)$ }
	\item[...]
	\item[n.] For any natural $n$ \\
		$\Delta_n = \{ v(m), v(m+1), v(m+2), ... \}$, \\ 
		where  $m = max\{K_k - N_{n-1,k-1} : k = 1,2,3,...,\lfloor n/2 \rfloor\}$
\end{enumerate}

\noindent where $K_k = \sum_{i=1}^{k}{\delta(v_{i})}$ and 
$ N_{n,k} = \sum_{i=n-k+1}^{n}{\delta(v_{i})}$.

\vspace{0.4cm}

\begin{defn}\label{def:family}
Denote with letter $\mathcal{A}(p)$ the family of all primary cobweb admissible sequences $P(p)$.
\end{defn}

\begin{observen} \label{obs:characterize}
The family $\mathcal{A}(p)$ is labelled-designated by the set of infinite paths $(v_0,v_1,v_2,...)$ of graph $G(p)$ from the root $v_0$ i.e.
\begin{center}
 $F\in \mathcal{A}(p) \Leftrightarrow (v_0,v_1,v_2,...)$ is a path of graph $G(p)$
\end{center}
where $F\equiv\{n_F\}_{n\geq 0}$ and $n_F = p^{\delta(v_n)}$.
\end{observen}

\noindent \textbf{P}\textbf{\footnotesize{ROOF}} \\
This is a conclusion on graph $G(p)$ (Definition \ref{def:tree}).\\
\textbf{The first steep}. If given any primary cobweb admissible sequence $F\equiv\{n_F\}_{n\geq 0}$, $n\in\{1,p,p^2,...\}$, then from the Definition \ref{def:admissible} of admissibility and the Definition \ref{def:tree} of tree $G(p)$ for any natural numbers $n,k$ the following is true
$$
		{n \choose k}_F = 
	\frac{p^{n_{B(F)}}\cdot...\cdot p^{(n-k+1)_{B(F)}}}
			{p^{1_{B(F)}} \cdot p^{2_{B(F)}}\cdot...\cdot p^{k_{B(F)}}} = 
	\frac{p^N}{p^K} =
	\frac{p^{\delta(v_n)}\cdot ... \cdot p^{\delta(v_{n-k+1})}}
			{p^{\delta(v_1}\cdot ... \cdot p^{\delta(v_k)}} 
	\in \mathbb{N} 
$$
where $s_{B(F)} = \delta(v_s)$ from Conclusion 1. In view of the Definition \ref{def:admissible} $N\geq K$  hence $(v_0,v_1,v_2,...)$ is a path of graph $G(p)$. 

\vspace{0.4cm}
\noindent \textbf{The second steep}. Take any given  path $(v_0,v_1,v_2,...)$ of the graph $G(p)$. Then by definition for any natural number $n,k$, $N_{n,k} \geq K_k$ where $N_{n,k} = \sum_{i=n-k+1}^{n}{\delta(v_{i})}$ and $K_k = \sum_{i=1}^{k}{\delta(v_{i})}$. Hence this path does encode the very  primary cobweb admissible sequence $P(p)$ $\blacksquare$

\begin{theoremn}[Cobweb Admissible Sequences Production Theorem]
The family of all cobweb admissible sequences is a product of families $\mathcal{A}(p_s)$ for $s=1,2,3,...$ i.e. for any cobweb admissible sequence $F$
$$
	F \in \times _{s=1}{\mathcal{A}(p_s)} 
$$
\end{theoremn}

\noindent \textbf{P}\textbf{\footnotesize{ROOF}} \\
\noindent This is the summarizing conclusion. Any cobweb admissible sequence $F$ is at the point product of primary cobweb admissible sequences $P(p)$ (Theorem \ref{th:1}) and the family of all primary cobweb admissible sequences $\mathcal{A}(p)$ is defined by primary cobweb admissible tree $G(p)$ (Observation \ref{obs:characterize}) $\blacksquare$

\vspace{0.4cm}
\noindent \textbf{Acknowledgements}

\vspace{0.4cm}
I would like to thank Professor A. Krzysztof Kwa\'sniewski for his supply of names for objects and operations and final improvements of this paper. Also discussions with Participants of Gian-Carlo Rota Polish Seminar \emph{http://ii.uwb.edu.pl/akk/sem/sem\_rota.htm} are appreciated.


\end{document}